\documentclass{amsart}
\usepackage{amssymb
,amsthm,dirtytalk
,amsmath
,amscd,mathrsfs
,mathtools,enumerate
}

\usepackage{hyperref}
\usepackage[all]{xy}
\usepackage[top=30truemm,bottom=30truemm,left=25truemm,right=25truemm]{geometry}

\newcommand{\GL}{\mathrm{GL}}

\newcommand{\Hom}{\mathrm{Hom}}

\newtheorem{thm}{Theorem}[section]

\newtheorem{lem}[thm]{Lemma}
\newtheorem{prop}[thm]{Proposition}

\theoremstyle{remark}
\newtheorem{rem}[thm]{Remark}

\theoremstyle{definition}

\numberwithin{equation}{section}

\makeatletter
\def\iddots{\mathinner{\mkern1mu\raise\p@
	\hbox{.}\mkern2mu\raise4\p@\hbox{.}\mkern2mu
	\raise7\p@\vbox{\kern7\p@\hbox{.}}\mkern1mu}}
\def\adots{\mathinner{\mkern2mu\raise\p@\hbox{.}
 \mkern2mu\raise4\p@\hbox{.}\mkern1mu
 \raise7\p@\vbox{\kern7\p@\hbox{.}}\mkern1mu}}
\makeatother

\allowdisplaybreaks
\title{Multiplicity one for the pair ($\GL_{n}(D),\GL_n(E)$) }
\author{Hengfei LU}
\address{Department of Mathematics, University of Vienna, Oskar-Morgenstern-Platz 1, Wien 1090, Austria}
\email{hengfei.lu@univie.ac.at}
\begin{document}

\begin{abstract} Let $F$ be a local field of characteristic zero. Let $D$ be a quaternion algebra over $F$. Let $E$ be a quadratic field extension of $F$. Let $\mu$ be a  character of $E^\times$. 
We study the distinction problem for the pair $(\GL_{n}(D),\GL_{n}(E))$ and we prove that any bi-$(\GL_n(E),\mu)$-equivariant tempered generalized function on $\GL_{n}(D)$ is invariant with respect to an anti-involution.
Then it implies that $\dim\Hom_{\GL_{n}(E)}(\pi,\mu)\leq 1$ by the generalized Gelfand-Kazhdan criterion. Thus we give a new proof to the fact that $(\GL_{2n}(F),\GL_{n}(E))$ is a Gelfand pair when $\mu$ is trivial and $D$  splits.
\end{abstract}
\maketitle
\tableofcontents

\section{Introduction}
Let $F$ be a local field of characteristic zero.  Let $E=F[\delta]$ be a quadratic field extension of $F$ with $\delta^2\in F^\times\setminus (F^\times)^2$. Let $Mat_{n,n}(E)$ denote the set of all $n\times n$ matrices over $E$. There is an embedding
\[\GL_{2n}(F)\hookrightarrow \GL_{2n}(E) \]
such that each element $g$ in the image of $\GL_{2n}(F)$ is of the form
\[\begin{pmatrix}
A& B\\ \bar{B}&\bar{A}
\end{pmatrix} \]
for $A,B\in Mat_{n,n}(E)$
 and the Galois action on $Mat_{n,n}(E)$ is given by
\[A\mapsto\bar{A}. \]
 (See \cite[Page 275]{guo1997unique}.)
 Let $\theta$ be the involution defined on $\GL_{2n}(F)$ given by
\[\theta(g)=\begin{pmatrix}
\delta\\&-\delta
\end{pmatrix} g\begin{pmatrix}
\delta\\&-\delta
\end{pmatrix}^{-1} \]
for $g\in\GL_{2n}(F)$.  Let $H_n$ be the fixed points of $\theta$ in $\GL_{2n}(F)$. Then $H_n\cong \GL_{n}(E)$ and it has elements of the form of block diagonal $\begin{pmatrix}
a\\&\bar{a}
\end{pmatrix}$, $a\in\GL_{n}(E)$. When $E=F\oplus F$, $H_n\cong \GL_n(F)\times \GL_n(F)$.  It is well known that the pair $(\GL_{2n}(F),\GL_n(F)\times\GL_n(F) )$ satisfies the Gelfand-Kazhdan criterion \cite[\S7]{dima2009duke} with respect to the inverse map; see \cite{jacquet1996linear} or \cite[Theorem 7.1.3]{dima2009duke}.
Following Jacquet-Rallis, if $F$ is nonarchimedean,  Guo \cite{guo1997unique} has proved that $(\GL_{2n}(F),\GL_n(E))$ satisfies the Gelfand-Kazhdan criterion with respect to the inverse map as well. Then Guo proved that $(\GL_{2n}(F),\GL_n(E))$ is a Gelfand pair, i.e.
\[\dim \Hom_{\GL_{n}(E)}(\pi,\mathbb{C})\leq1 \]
for any irreducible smooth representation $\pi$ of $\GL_{2n}(F)$.
Broussons-Matringe recently \cite{matringe19} verified a very general version that  the inner form  of the pair $(\GL_{2n}(F),\GL_n(E))$ is a Gelfand-pair. 
 The case when $F=\mathbb{R}$ has been proved by S. Carmeli in \cite{carmeli2015stability}. This paper studies the twisted version.

The main result in this paper is the following:
\begin{thm}\label{thm1.1} 
	 Given any irreducible admissible smooth representation $\pi$ of $\GL_{2n}(F)$, we have
	 \[\dim \Hom_{\GL_{n}(E)}(\pi,\mu\circ\det)\leq1 \]
	 for any  character $\mu $ of $E^\times$.
\end{thm}

The key point is that  the symmetric pair $(\GL_{2n}(F),\GL_n(E))$  satisfies the generalized Gelfand-Kazhdan criterion with respect to not only the inverse map, which implies that it is a Gelfand pair, but also the transposition and then \cite[Theorem 2.3]{sunzhu2011} implies that Theorem \ref{thm1.1} holds. 
\begin{thm}\label{thm:A}
	Let $f$ be  a tempered generalized function on $\GL_{2n}(F)$. Let $\mu=\mu\circ\det$ be a  character of $\GL_{n}(E)$.
	 If for every $h\in \GL_{n}(E)$,
	\[f(hx)=f(xh)=\mu(h)f(x) \]
for $x\in\GL_{2n}(F)$,	as generalized functions on $\GL_{2n}(F)$, then 
\[f(x)=f(x^t). \]
\end{thm}
Here and as usual, a superscript \say{t} indicates the transpose of a matrix. 
\begin{rem}
	Theorem \ref{thm1.1} generalizes a result of Waldspurger   \cite{waldspurger} for $n=1$. Using a new anti-involution, we give a new proof
	to the result that $(\GL_{2n}(F),\GL_n(E))$ is a Gelfand pair when $\mu$ is trivial. There is a precise conjecture, in terms of the Langlands parameter and the $\epsilon$-factor, to predict when the multiplicity is $1$. (See \cite[Conjecture 1]{prasadramin2008crelle} or \cite[Conjecture 4]{dipendraJNT}.) There are a lot of progress for this conjecture recently (see \cite{feigonkimball,xuesuzuki,vincent20,xuehang}) which will not be addressed in this paper.
	In fact, \cite[Conjecture 1]{prasadramin2008crelle} is  very general which is   related to the pair $(G',H')$ where $G'$ (resp. $H'$) is an inner form of $\GL_{2n}(F)$ (resp. $\GL_n(E)$). 
	This paper focuses on the pair $(\GL_{n}(D),\GL_n(E))$, where $D$ is a quaternion algebra of $F$. Thus it gives a partial answer to the question in \cite[Remark 11.4]{prasadramin2008crelle}
\end{rem}

 The analogue for the pair $(\GL_{2n}(F),\GL_n(F)\times\GL_n(F) )$ has been proved for \say{good} characters by Chen-Sun (see \cite{sun2020} for more details). We will use a similar idea appearing in \cite{sun2020} to prove Theorem \ref{thm:A}.

By linearization, Theorem \ref{thm:A} is reduced to the following theorem.
\begin{thm}\label{tran:M} 
	Let $f$ be a tempered generalized function on $Mat_{n,n}(E)$ 
	 such that for any $h$ in $\GL_n(E)$,
	\[f(hx\bar{h}^{-1})=f(x) \]
holds	for any $x\in Mat_{n,n}(E)$. Then $f(x)=f(\bar{x}^t)$.
\end{thm} 
Here is a brief introduction to the proof of Theorem \ref{tran:M}. It is enough to show that if a tempered generalized function $f$ satisfies $f(hx\bar{h}^{-1})=f(x)$ and $f(x)=-f(\bar{x}^t)$, then $f=0$. We will use Aizenbud-Gourevitch's result (see Theorem \ref{thm2.1}) to reduce the problem to the nilpotent cone. Inside the nilpotent cone,
we will pick up an $\mathfrak{sl}_2(F)$-triple (see \eqref{sl}) and use the result of Chen-Sun (see Lemma \ref{key:lem})   to prove the vanishing theorem (see Theorem \ref{vanishing:I}).

Let $D$ be a quaternion algebra over $F$. If $D$ splits, then
$\GL_n(D)=\GL_{2n}(F)$.
	All the techniques in this paper work for the pair $(\GL_{n}(D),\GL_{n}(E))$ as well when $D$ is the unique quaternion division algebra over $F$. Fix an element $\epsilon$ in $F^\times$ which does not belong to $N_{E/F}E^\times$. Replace the transpose by another anti-involution
	\[\sigma':g\mapsto\begin{pmatrix}
	\epsilon&\\&1
	\end{pmatrix}g^t\begin{pmatrix}
	\epsilon\\&1
	\end{pmatrix}^{-1} \]
	for all $g\in \GL_{n}(D)$ which consists of matrices of the form $\begin{pmatrix}
	A&\epsilon B\\ \bar{B}&\bar{A}
	\end{pmatrix}$. (See \cite{guo1997unique}.) Then the proof will follow  in a similar way (see \S\ref{sect:inner}). However, this method does not work for other inner form pairs ($G',H'$) where $G'$ (resp. $H'$)
	is an inner form of $\GL_{2n}(F)$ (resp. $\GL_n(E)$) corresponding to a central simple algebra (containing $E$) over $F$ and $H'$ corresponds to the centralizer of $E$ inside the central simple algebra. One reason is that there is lack of an involution of $G'$ which takes any irreducible representation $\pi$ of $G'$ to its contragredient $\pi^\vee$. Another reason is that there is lack of  a proper anti-involution of $G'$ keeping $H'$
	which intertwines the Fourier transform, while the anti-involution $\sigma:g\mapsto g^t$ (resp. $\sigma'$) of $\GL_{2n}(F)$ (resp. $\GL_{n}(D)$) plays a vital role in the Fourier transform and thus in the proof of Theorem \ref{thm:A} (resp. Theorem \ref{G'}).

The paper is organized as follows. In \S2, we introduce some notation from the algebraic geometry. Then we will use Chen-Sun's method to prove Theorem \ref{tran:M} in \S3. The proof of Theorem \ref{thm:A} will be given 
in \S4 and then Theorem \ref{thm1.1} will follows easily. We will study the multiplicity one problem for the pair $(\GL_{n}(D),\GL_n(E))$ in the last section.

\section{Preliminaries and notation}

Let $X$ be an $\ell$-space (i.e. locally compact totally disconnected topological spaces) if $F$ is nonarchimedean, or a Nash manifold if $F=\mathbb{R}$ (see \cite[\S2.3]{dima2009duke}).  Let $\mathscr{C}(X)$ denotes the tempered generalized functions on $X$ (see \cite{dima2009duke} for the precise definition). Let a reductive group $G(F)$
act on an affine variety $X$. Given a character $\chi$ of $G(F)$, denote by $\mathscr{C}(X)^{G(F),\chi}$
the $(G(F),\chi)$-equivariant tempered generalized functions on $X$. Let $x\in X$ such that its orbit $G(F)x$ is closed in $X$. We denote the normal bundle by $N_{G(F)x,x}^X$. Let $$G_x:=\{g\in G(F)|gx=x \}$$ be the stalizer subgroup of $x$.
\begin{thm}\cite[Theorem 3.1.1]{dima2009duke}\label{thm:duke}
	Let $G(F)$ act on a smooth affine variety $X$. Let $\chi$ be a character of $G(F)$. Suppose that for any closed orbit $G(F)x$ in $X$, we have
	\[\mathscr{C}(N^X_{G(F)x,x})^{G_x,\chi}=0. \]
	Then \[\mathscr{C}(X)^{G(F),\chi}=0. \]
\end{thm}
 If $V$ is a finite
dimensional representation of $G(F)$, then we denote the nilpotent cone in $V$ by
\[S_0:=\{ x\in V|\overline{G(F)x}\owns 0 \}. \] 
Let $Q_G(V):=V/V^G$ which can be regarded as a subset in $V$ (see \cite[Notation 2.3.10]{dima2009duke}) and $R_G(V):=Q_G(V)\setminus S_0$. There is a stronger version of Theorem \ref{thm:duke}.
\begin{thm}\cite[Corollary 3.2.2]{dima2009duke}\label{thm2.1}
	Let $G(F)$ act on a smooth affine variety $X$. Let $K\subset G(F)$ be an open subgroup and let $\chi$ be a character of $K$. Suppose that for any closed orbit $G(F)x$ such that
	\[\mathscr{C}(R_{G_x}(N_{G(F)x,x}^X))^{K_x,\chi}=0 \]
	we have \[\mathscr{C} (Q_{G_x}(N^X_{G(F)x,x}))^{K_x,\chi}=0. \]
	Then $\mathscr{C}(X)^{K,\chi}=0$.
\end{thm}

\section{A vanishing result of generalized functions}
In this section, we shall prove Theorem \ref{tran:M}.
Define
$$L_n:=Mat_{n,n}(E)=\Big\{\begin{pmatrix}
0&x\\ \bar{x}&0
\end{pmatrix}:x\in Mat_{n,n}(E) \Big\}\subset \mathfrak{gl}_{2n}(F).$$
Denote by 
\[\mathcal{N}_n:=\{\begin{pmatrix}
	0&x \\
	\bar{x}&0
\end{pmatrix}\in L_n|x \bar{x} \mbox{ is a nilpotent matrix in }Mat_{n,n}(E) \} \]
the nilpotent cone in $L_n$. Let $\mathscr{C}_{\mathcal{N}_n}(L_n)$ denote the tempered generalized functions supported on $\mathcal{N}_n$.
 Denote $\tilde{H}_n:=H_n\rtimes\langle\sigma\rangle$ where $H_n=\GL_n(E)$ and $\sigma$ acts on $H_n$ by the involution 
$$\begin{pmatrix}
a\\&\bar{a}
\end{pmatrix}\mapsto\begin{pmatrix}
(a^{-1})^t\\&(\bar{a}^{-1})^t
\end{pmatrix}.$$
 The group $\tilde{H}_n$ acts on $L_n$ by
\[\begin{pmatrix}
a\\&\bar{a}
\end{pmatrix}\cdot \begin{pmatrix}
0&x\\\bar{x}&0
\end{pmatrix}=\begin{pmatrix}0&ax\bar{a}^{-1}\\ \bar{a}\bar{x}a^{-1}&0\end{pmatrix} \]
and $$\sigma\cdot \begin{pmatrix}
	0&x\\\bar{x}&0
\end{pmatrix}=\begin{pmatrix}
0&\bar{x}^t\\ x^t &0
\end{pmatrix}$$
 for $x\in Mat_{n,n}(E)$. Let $\chi$ be the sign character of $\tilde{H}_n$ i.e.,
$\chi|_{H_n}$ is trivial  and
\[\chi(\sigma)=-1. \]
\begin{thm}\label{vanishing:I}
	We have $\mathscr{C}(L_n)^{\tilde{H}_n,\chi}=0$.
\end{thm}
The rest of this section is devoted to proving Theorem \ref{vanishing:I}. Indeed Theorem \ref{vanishing:I} is  a reformulation of  Theorem \ref{tran:M}.

Define a non-degenerate symmetric $F$-bilinear form on $\mathfrak{gl}_{2n}(F)$ by
\[\langle z, w\rangle_{\mathfrak{gl}_{2n}(F)}:=\mbox{the trace of }zw\mbox{ as a }F\mbox{-linear operator}. \]
Note that the restriction of this bilinear form on $L_n$ is still non-degenerate. Fix a non-trivial unitary character $\psi$ of $F$. Denote by
\[\mathfrak{F}:\mathscr{C}(L_n)\longrightarrow \mathscr{C}(L_n) \]
the Fourier transform which is normalized such that for every Schwartz function $\varphi$ on $L_n$,
\[\mathfrak{F}(\varphi)(z)=\int_{L_n}\varphi(w)\psi(\langle z,w \rangle_{\mathfrak{gl}_{2n}(F)})d w  \]
for $z\in L_n$, where $d w$ is the self-dual Haar measure on $L_n$. More precisely
\[\langle \begin{pmatrix}
	0&x\\ \bar{x}&0
\end{pmatrix},\begin{pmatrix}
0&x'\\ \bar{x}'&0
\end{pmatrix}\rangle_{\mathfrak{gl}_{2n}(F)}=tr(x\bar{x}'+x'\bar{x})=N_{E/F}tr(x\bar{x}') \]
for $x,x'\in Mat_{n,n}(E)$ and $x\bar{x}'+x'\bar{x}\in Mat_{n,n}(F)$.
 It is clear that the Fourier transform $\mathfrak{F}$ intertwines the action of $\tilde{H}_n$. Thus we have the following lemma.
\begin{lem}
	The Fourier transform $\mathfrak{F}$ preserves the space $\mathscr{C}(L_n)^{\tilde{H}_n,\chi}$.\label{lem:intertwin}
\end{lem}

\subsection{Reduction within the null cone} We will use $(x,\bar{x})$ to denote the matrix $\begin{pmatrix}
	0&x\\ \bar{x}&0
\end{pmatrix}$ in $L_n$.
Recall that 
\[\mathcal{N}_n=\{(x,\bar{x})\in L_n|x \bar{x} \mbox{ is a nilpotent matrix in }Mat_{n,n}(E) \} \]
the nilpotent cone in $L_n$.
\begin{lem}\cite[Lemma 2.3]{guo1997unique}
	Let $(x,\bar{x})\in L_n$ be nilpotent. There exists an element $\begin{pmatrix}
	a\\&\bar{a}
	\end{pmatrix}$ in $H_n$ such that the twisted conjugate $ax\bar{a}^{-1}$ of $x$ is in its Jordan normal form.\label{guotwisted}
\end{lem}
 Let $\mathcal{O}$ be an $H_n$-orbit in $\mathcal{N}_n$. Recall that every $\mathbf{e}\in\mathcal{O}$ can be extended to a  $\mathfrak{sl}_2$-triple $\{\mathbf{h},\mathbf{e},\mathbf{f} \}$ (see \cite[Proposition 4]{kostant}) in the sense that 
\begin{equation}\label{sl}
[\mathbf{h},\mathbf{e}]=2\mathbf{e},~[\mathbf{h},\mathbf{f}]=-2\mathbf{f}\mbox{ and }[\mathbf{e},\mathbf{f}]=\mathbf{h} \end{equation}
where $\mathbf{f}\in\mathcal{N}_n$ and $\mathbf{h}\in\mathfrak{h}$, where $\mathfrak{h}=\mathfrak{gl}_{n}(E)$ is the Lie algebra of $H_n$. Let $L_n^{\mathbf{f}}$ denote the elements in $L_n$ annihilated by $\mathbf{f}$ under the adjoint action of the  $\mathfrak{sl}_2$-triple $\{\mathbf{h},\mathbf{e},\mathbf{f} \}$ on $L_n\subset\mathfrak{gl}_{2n}(F)$. Then 
\[L_n=[\mathfrak{h},\mathbf{e}]\oplus L_n^{\mathbf{f}}. \]
Moreover, we may assume that $\mathbf{e}=(x,\bar{x})$ such that $x\in Mat_{n,n}(F)$ (so that $x$ is nilpotent) and $\mathbf{h}\in\mathfrak{gl}_n(F)$ due to Lemma \ref{guotwisted} which states that the orbit $\mathcal{O}$ always contains an element in $Mat_{n,n}(F)$.

 Define $\mu_1(g):=\mu(\det g^{-1}\bar{g})$ for $g\in\GL_n(E)$.
 Following \cite[Proposition 3.9]{sun2020},
we shall prove the following proposition in this subsection.
\begin{prop}\label{fourier} 
	Let $f$ be a $(H_{n},\mu_1)$-equivariant tempered generalized function on $L_{n}$ such that 
	 $f$  and its Fourier transforms $\mathfrak{F}(f)$  are both supported on  $\mathcal{N}_n$.  Then $f=0$.
\end{prop}

Denote by $\mathscr{C}_\mathcal{O}(L_{n})$ the space of tempered generalized functions on $L_{n}\setminus (\partial \mathcal{O})$ with support in $\mathcal{O}$, where $\partial \mathcal{O}$ is the complement of $\mathcal{O}$ in its closure in $L_{n}$. (See \cite[Notation 2.5.3]{dima2009duke}.) We will  use similar notation without further explaination.

Let $F^\times$ act on $\mathscr{C}(L_{n})$ by
\[(t\cdot f)(x)=f(t^{-1}x) \]
for $t\in F^\times$, $x\in L_{n}$ and $f\in\mathscr{C}(L_{n})$. The orbit $\mathcal{O}$ is invariant under dilation and so $F^\times$ acts on $\mathscr{C}_{\mathcal{O}}(L_{n})^{H_{n},\mu_1}$ as well.
\begin{lem}\cite[Lemma 3.13]{sun2020}
	Let $\eta:F^\times\rightarrow\mathbb{C}^\times$ be an eigenvalue for the action of $F^\times$ on $\mathscr{C}_{\mathcal{O}}(L_{n})^{H_{n},\mu_1}$. Then
	$\eta^2=|-|^{tr(2-\mathbf{h})|_{L_{n}^\mathbf{f}}}\cdot\kappa$, where $|-|$ denotes the absolute value of $F^\times$, $\kappa=|-|^r$ for some non-negative integer $r$ when $F=\mathbb{R}$ and $\kappa$ is the trivial character when $F$ is non-archimedean. 
	\label{key:lem}
\end{lem}
Note that the restrictions on the eigenvalues established in the lemma do not depend on the character $\mu_1$.

Let $Q$ be a quadratic form on $L_{n}$ defined by
\[Q(x)=2tr(x\bar{x}) \]
for $x\in L_{n}$. Denote by $Z(Q)$ the zero locus of $Q$  in $L_n$. Then $\mathcal{N}_{n}\subset Z(Q)\subset L_{n}$. Recall the following homogeneity result on tempered generalized functions. (See \cite[Theorem 2.1.10]{aizenbud2013partial}.)
\begin{thm}\label{duke}\cite[Proposition 3.14]{sun2020}
		Let $I$ be a non-zero subspace of $\mathscr{C}_{Z(Q)}(L_{n})$ such that for every $f\in I$, one has that $\mathfrak{F}(f)\in I$ and $(\psi\circ Q)\cdot f\in I$ for all unitary character $\psi$ of $F$. Then $I$ is a completely reducible $F^\times$-subrepresentation of ${\mathscr{C}}(L_{n})$, and it has an eigenvalue of the form 
	\[\begin{cases}
	|-|^{\frac{1}{2}\dim_F L_n},&\mbox{ if }F\mbox{ is non-archimedean},\\
	|-|^{\frac{1}{2}\dim_F L_n-c},&\mbox{ if }F=\mathbb{R},
	\end{cases}\]
	for a non-negative integer $c$.
\end{thm}

Now we are prepared to prove Proposition \ref{fourier}. The basic idea is due to Chen-Sun  \cite{sun2020}.
\begin{proof}[Proof of Proposition \ref{fourier}]
	Denote by $I$ the space of all tempered generalized functions $f$ on $L_{n}$ with the properties in Proposition \ref{fourier}. Assume by contradiction that $I$ is nonzero. We separate the proof into two cases.
	\begin{enumerate}[(i)]
		\item
	Suppose that $F$ is nonarchimedean. Then by Lemma \ref{key:lem} and Theorem \ref{duke}, one has
	\[|-|^{\dim_F L_{n}}=\eta^2=|-|^{tr(2-\mathbf{h})|_{L_{n}^\mathbf{f}}}. \]
	Thus $2n^2=tr(2-\mathbf{h})|_{L_n^\mathbf{f}}$. Thanks to Lemma \ref{guotwisted}, we may  assume that $\mathbf{e}=(x,x)$ where $x$ is an upper triangular nilpotent matrix in $Mat_{n,n}(F)$. Similarly  $\mathbf{f}=(x',x')$  where $x'$ is a lower triangular nilpotent matrix in $Mat_{n,n}(F)$
	and $\mathbf{h}=[\mathbf{e},\mathbf{f}]$.
	 Consider $Mat_{n,n}(F)$ as a representation of this $\mathfrak{sl}_2(F)$-triple. Decompose it into irreducible representations
	\[Mat_{n,n}(F)=\oplus_{i=1}^d V_i .\]
	Let $\lambda_i$ be the highest weight of $V_i$. Then
	\[2n^2=tr(2-\mathbf{h})|_{L_n^\mathbf{f}}=2(2d+\sum_{i=1}^d\lambda_i) \]
	and $n^2=\sum_i(\lambda_i+1)=d+\sum_i\lambda_i$ (see \cite[Lemma 7.6.6]{dima2009duke}). Therefore we have
	\[d+\sum_i\lambda_i=2d+\sum_i\lambda_i \]
which is impossible. 
\item
Now we suppose that $F=\mathbb{R}$. Then Lemma \ref{key:lem} and Theorem \ref{duke} imply that
\begin{equation}\label{key:arch}
\dim_F L_n-c=tr(2-\mathbf{h})|_{L_n^\mathbf{f}}+r 
\end{equation}
for some non-negative integers $c$ and $r$ . However,  
\[\dim_F L_n-tr(2-\mathbf{h})|_{L_n^\mathbf{f}}=2\sum_{i=1}^d(\lambda_i+1)-2( 2d+\sum_{i=1}^d\lambda_i)=-2d<0\]
and so \eqref{key:arch} must be false.
\end{enumerate}
This finishes the proof.
\end{proof}

\subsection{Proof of Theorem \ref{vanishing:I}} We shall use the generalized Harish-Chandra descent developed by Aizenbud-Gourevitch in \cite{dima2009duke} to prove Theorem \ref{vanishing:I} in this subsection.

Let $(G,H,\theta)$ be a symmetric pair. Consider the action of $H\times H$ on $G$ by left and right translations and the action of $H$ on $L_n$ by conjugation. Let $g\in G(F)$ such that $H(F)gH(F)$ is closed in $G(F)$. Let $x=g\theta(g^{-1})$. Then
$(G_x,H_x,\theta|_{G_x})$ is called a descendant of $(G,H,\theta)$. (See \cite[\S7.2]{dima2009duke}.)
\begin{lem}\cite[Theorem 6.15]{carmeli2015stability}\label{descendant}
	Every descendant of the pair $(\GL_{2n},R_{E/F}\GL_{n})$ is a product of pairs of the form $(R_{K_1/F}\GL_{r},R_{K_2/F}\GL_r)$, $(\GL_r\times\GL_r,\triangle\GL_r)$ and  $(\GL_{2r},R_{E/F}\GL_r)$ for some $r<n$, where $K_2$ is a finite field extension of $F$ and $K_1$ is a quadratic field extension of $K_2$.
\end{lem}
Here the direct summand of the normal space at the closed orbit associated to the descendant $(R_{K_1/F}\GL_r,R_{K_2/F}\GL_r)$ corresponds to the linear space
$\delta'Mat_{r,r}(K_2)$ with the congugation action of $\GL_r(K_2)$ where $\delta'\in K_1$ and $tr_{K_1/K_2}(\delta')=0$, i.e. 
\[g\cdot \delta'x=\delta'gxg^{-1} \]
for $g\in\GL_r(K_2)$ and $x\in Mat_{r,r}(K_2)$.

Now we can give a proof of Theorem \ref{vanishing:I}.

\begin{proof}[Proof of Theorem \ref{vanishing:I}] 
	Define $\widetilde{\GL_r(F)}:=\GL_r(F)\rtimes\langle\sigma\rangle$ with $\sigma(g)=(g^t)^{-1}$ for $g\in\GL_r(F)$ and $\sigma(y)=y^t$ for $y\in Mat_{r,r}(F)$. Here $\GL_r(F)$ acts on $Mat_{r,r}(F)$
	by conjugation. Let $\chi$ be the sign character of $\widetilde{\GL_r(F)}$, i.e. $\chi(g)=1$ for $g\in\GL_r(F)$ and $\chi(\sigma)=-1$.
	Thanks to \cite[Theorem D(i)]{sun2020}, one has $$\mathscr{C}(Mat_{r,r}(F))^{\widetilde{\GL_r(F)},\chi}=0.$$
	It implies that the descendants $(R_{K_1/F}\GL_r,R_{K_2/F}\GL_r )$ and $(\GL_r(F)\times\GL_r(F),\triangle\GL_r(F))$ do not contribute to $\mathscr{C}(L_n)^{\tilde{H}_n,\chi}$.
	Applying Theorem \ref{thm2.1} and Lemma \ref{descendant},
	it suffices to show that ${\mathscr{C}}(R(L_{r}))^{\tilde{H}_{r},\chi}=0\Longrightarrow \mathscr{C}(Q(L_r))^{\tilde{H}_r,\chi}=0$ for all $r$. Then
 Theorem \ref{vanishing:I} follows from Lemma \ref{lem:intertwin} and Proposition \ref{fourier} directly. More precisely, if $f\in \mathscr{C}(L_r)^{\tilde{H}_r,\chi}$ supported on the nilpotent cone $\mathcal{N}_r$, then its Fourier transform $\mathfrak{F}(f)$ is also supported on $\mathcal{N}_r$. Thus $f=0$ due to Proposition \ref{fourier}.
\end{proof}
\begin{proof}[Proof of Theorem \ref{tran:M}]
	Note that $\sigma$ acts on $Mat_{n,n}(E)$ with order $2$. Then $\sigma\cdot(f-\sigma\cdot f)=\sigma\cdot f-f$ for any $f\in\mathscr{C}(Mat_{n,n}(E))^{H_n}$. Then $f-\sigma\cdot f=0$ due to Theorem \ref{vanishing:I}, i.e.
	$f(x)=f(\bar{x}^t)$ for any $x\in Mat_{n,n}(E)$ and $f\in\mathscr{C}(Mat_{n,n}(E))^{H_n}$.
\end{proof}

\section{Proof of Theorem \ref{thm:A}} 
Let $\mathcal{H}_{n}:=H_{n}\times H_{n}$ be a  reductive group.  Define
\[\tilde{\mathcal{H}}_{n}:=\mathcal{H}_{n}\rtimes\langle\sigma\rangle \]
where $\sigma$
acts on $\mathcal{H}_{n}$ by the involution $(h_1,h_2)\mapsto ((h_2^{-1})^t,(h_1^{-1})^t)$. Let $\tilde{\mathcal{H}}_{n}$ act on $\GL_{2n}(F)$ by
\[(h_1,h_2)\cdot g=h_1gh_2^{-1} \]
and $\sigma\cdot g=g^t$ for $h_i\in H_{n}$ and $g\in\GL_{2n}(F)$.
Let $\mu$ be a  character of $H_n$  and let $\mu\otimes{\mu}^{-1}$ be a character of $\mathcal{H}_n$. Let $\tilde{\mu}$ be the character of $\tilde{\mathcal{H}}_n$ twisted by the sign character, i.e.
$\tilde{\mu}|_{\mathcal{H}_n}=\mu\otimes{\mu}^{-1}$ and $\tilde{\mu}(\sigma)=-1$.

This section is devoted to  the proof of the following theorem.
\begin{thm}
	\label{thm:dist} We have
	\[{\mathscr{C}}(\GL_{2n}(F))^{\tilde{\mathcal{H}}_{n},\tilde{\mu}}=0. \]
\end{thm}
Then Theorem \ref{thm:A} will follow from Theorem \ref{thm:dist} immediately.

\begin{lem}\label{Jac:lem} \cite[Proposition 1.2]{guo1996cjm}
	 The double cosets $H_nxH_n$, where
	$$x\theta(x^{-1})=\begin{pmatrix}
	A&0&0&B_A&0&0\\0&-\mathbf{1}_{p}&0&0&0&0\\ 0&0&\mathbf{1}_q&0&0&0\\
	\bar{B}_A&0&0&A&0&0\\
	0&0&0&0&-\mathbf{1}_p&0\\
	0&0&0&0&0&\mathbf{1}_q
	\end{pmatrix},$$
	exhaust all closed orbits
	 in $H_n\backslash \GL_{2n}(F)/H_n$, where $A\in Mat_{\nu,\nu}(F)$ is semisimple without eigenvalues $\pm1$,
	  $\nu+p+q=n$,  $\mathbf{1}_p,\mathbf{1}_q,\mathbf{1}_\nu$ are identity matrices and $B_A\in Mat_{\nu,\nu}(E)$ satisfies
	$A^2-\mathbf{1}_{\nu}=B_A\bar{B}_A$  and $AB_A=B_A A$. 
\end{lem}
\begin{proof}
	See \cite[Theorem 4.13]{carmeli2015stability}.
\end{proof}

Thanks to Theorem \ref{thm2.1}, if 
\[\mathscr{C}(R(N_{\mathcal{O},x}^{\GL_{2n}(F)}))^{\tilde{\mathcal{H}}_{n,x},\tilde{\mu}}=0\Longrightarrow \mathscr{C}(Q(N_{\mathcal{O},x}^{\GL_{2n}(F)}))^{\tilde{\mathcal{H}}_{n,x},\tilde{\mu}}=0 \]
for any $\tilde{\mathcal{H}}_{n}$-closed orbit $\mathcal{O}=\tilde{\mathcal{H}}_{n}x$, where $\tilde{\mathcal{H}}_{n,x}$ is the stabilizer of $x$, then Theorem \ref{thm:dist} holds.

We will divide the proof of Theorem \ref{thm:dist} into two cases: $\nu=0$ and $\nu\neq0$.
\begin{lem}\label{nu=0} Suppose that $\nu=0$ and $x=\begin{pmatrix}
&&	\mathbf{1}_p\\&\mathbf{1}_q\\\mathbf{1}_p\\&&&\mathbf{1}_q
	\end{pmatrix}$ with $p+q=n$.
	We have 
	$${\mathscr{C}}(R(N^{\GL_{2n}(F)}_{\tilde{\mathcal{H}}_nx,x }))^{\tilde{\mathcal{H}}_{n,x},\tilde{\mu}}=0\Longrightarrow {\mathscr{C}}(Q(N^{\GL_{2n}(F)}_{\tilde{\mathcal{H}}_nx,x }))^{\tilde{\mathcal{H}}_{n,x},\tilde{\mu}}=0  $$
	where $\tilde{\mathcal{H}}_{n,x}=\{h\in\tilde{\mathcal{H}}_{n}|h\cdot x=x \}$ is the stabilizer of $x$ in $\tilde{\mathcal{H}}_n$.
\end{lem}
\begin{proof}
	By easy computation, we have $\tilde{\mathcal{H}}_{n,x}\cong (\GL_p(E)\times\GL_q(E))\rtimes\langle \sigma\rangle$. Moreover, $\tilde{\mu}|_{\GL_q(E)\rtimes\langle\sigma\rangle }=\chi$ is the sign character and $\tilde{\mu}|_{\GL_p(E)}=\mu_1$.
	 The normal space (see \cite[Lemma 4.3]{sun2020}) is given by
	\[N_{\tilde{\mathcal{H}}_{n}x,x}^{\GL_{2n}(F)}= \frac{\mathfrak{gl}_{2n}(F)}{\mathfrak{h}+Ad_{x}\mathfrak{h}{} }\cong Mat_{p,p}(E)\oplus Mat_{q,q}(E). \]
The action of $\tilde{\mathcal{H}}_{n,x}$ on $N_{\tilde{\mathcal{H}}_nx,x}^{\GL_{2n}(F)}$ is given by
\[(g_1,g_2)\cdot (x,y)=(\bar{g}_1xg_1^{-1},g_2y\bar{g}_2^{-1}) \]
and $\sigma\cdot(x,y)=(x^t,\bar{y}^t)$ for $g_1\in\GL_p(E),g_2\in \GL_q(E),x\in Mat_{p,p}(E)$ and
$y\in Mat_{q,q}(E)$. 
Thanks to \cite[Proposition 2.5.8]{dima2009duke} and Theorem \ref{vanishing:I}, it is enough to show that
any $(H_p,\mu_1)$-equivariant tempered generalized functions on $Mat_{p,p}(E)$ supported on the nilpotent cone $\mathcal{N}_p$ is invariant under the transposition, 
i.e. $f(\bar{g}xg^{-1})=\mu_1(g)f(x)\Longrightarrow f(x^t)=f(x)$ for any $x\in\mathcal{N}_p$ where $g\in H_p$. Consider $\mathscr{C}(Mat_{p,p}(E))^{H_p,\mu_1}$ and its Fourier transform
\[\mathfrak{F}(\varphi)(z)=\int_{Mat_{p,p}(E)}\varphi(y)\psi\circ N_{E/F}(tr(z\bar{y}))dy  \]
 which intertwines the action of $\tilde{H}_p$. Moreover, the Fourier transform preserves the space
 \[\mathscr{C}_{\mathcal{N}_p}(Mat_{p,p}(E))^{H_p,\mu_1. } \]
We may choose an $\mathfrak{sl}_2$-triple $\{\mathbf{h},\mathbf{e},\mathbf{f} \}$ lying in $Mat_{p,p}(F)$ due to Lemma \ref{guotwisted}. Thanks to Proposition \ref{fourier}, we obtain
that any $(H_p,\mu_1)$-equivariant tempered generalized function $f$ on $Mat_{p,p}(E)$ such that both $f$ and its Fourier transform $\mathfrak{F}(f)$ are supported on $\mathcal{N}_p$ must be zero.
Then
 we have finished the proof.
\end{proof}

Now we can give a proof of Theorem \ref{thm:dist}.
\begin{proof}[Proof of Theorem \ref{thm:dist} ] Applying Theorem \ref{thm2.1}, we only need to prove that there does not exist any $(\tilde{\mathcal{H}}_{n},\tilde{\mu})$-equivariant tempered generalized function supported on the nilpotent cone of the normal space of $\tilde{\mathcal{H}}_{n}$-closed orbits.
	Thanks to Lemma \ref{nu=0}, we have proved  Theorem \ref{thm:dist} if $\nu=0$. Thus applying Lemma \ref{Jac:lem}, it is reduced to prove that 
	\begin{equation}\label{substep}
	\mathscr{C}(R(N_{\tilde{\mathcal{H}}_nx,x}^{\GL_{n}(F)}))^{\tilde{\mathcal{H}}_{n,x},\tilde{\mu}}=0\Longrightarrow
	{\mathscr{C}}(Q(N^{\GL_{2n}(F)}_{\tilde{\mathcal{H}}_nx,x }))^{\tilde{\mathcal{H}}_{n,x},\tilde{\mu}}=0 \end{equation}
	for the closed orbit $\tilde{\mathcal{H}}_{n}x$, where $x$ satisfies
	\[x\theta(x^{-1})=\begin{pmatrix}
	A&&B_A\\&\mathbf{1}_{n-\nu}&\\ \bar{B}_A&&A\\&&&\mathbf{1}_{n-\nu}
	\end{pmatrix}\in\GL_{2n}(F), \]
	$A^2-\mathbf{1}_\nu=B_A\bar{B}_A,AB_A=B_AA$ and 
	 $A$ is a semisimple element in $Mat_{\nu,\nu}(F)$ for $\nu=1,2,\cdots n,$ without eigenvalues $\pm1$.
	Futhermore, we may assume that $A$ is a scalar matrix and $A^2\neq\mathbf{1}_\nu$. Then 
	\[\tilde{\mathcal{H}}_{n,x}\cong(\GL_{\nu}(F)\times \GL_{n-\nu}(E))\rtimes\langle\sigma\rangle, \]
$\tilde{\mu}|_{\tilde{\mathcal{H}}_{n,x}}=\chi$ is the sign character	and $N_{\tilde{\mathcal{H}}_{n}x,x}^{\GL_{2n}(F)}\cong Mat_{\nu,\nu}(F)\oplus L_{n-\nu}$, where $\GL_{\nu}(F)$ acts on $Mat_{\nu,\nu}(F)$ by inner conjugation and $\sigma$ acts on $Mat_{\nu,\nu}(F)$ by the matrix transpose. (See \cite[Proposition 7.2.1]{dima2009duke}.) Therefore \eqref{substep} follows from 
	\[\mathscr{C}(Mat_{\nu,\nu}(F))^{\GL_\nu(F)\rtimes\langle\sigma\rangle,\chi }=0=\mathscr{C}(L_{n-\nu})^{\tilde{\mathcal{H}}_{n-\nu},\chi}. \]
 (see \cite[Theorem D(i)]{sun2020} and Theorem \ref{vanishing:I}).
 	This finishes the proof.
\end{proof}

Finally, we can give a proof of Theorem \ref{thm1.1}.
\begin{proof}[Proof of Theorem \ref{thm1.1}]
 Theorem \ref{thm:A} and \cite[Theorem 2.3]{sunzhu2011} imply that
\[\dim \Hom_{\GL_{n}(E)}(\pi,\mu)\cdot\dim\Hom_{\GL_{n}(E)}(\pi^\vee,{\mu}^{-1})\leq1 \]
for any irreducible smooth admissible representation $\pi$ of $\GL_{2n}(F)$ and any  character $\mu$ of $E^\times$.
 Note that
\[\dim\Hom_{\GL_{n}(E)}(\pi^\vee,\mu^{-1})=\dim \Hom_{\GL_{n}(E)}(\pi,\mu) \]
due to the fact that the Chevalley involution  takes $\pi$ to its contragredient (see \cite{ka} when $F$ is a p-adic field and \cite{adams} when $F=\mathbb{R}$). 
Thus $\dim\Hom_{\GL_{n}(E)}(\pi,\mu)\leq 1 $, i.e. Theorem \ref{thm1.1} holds.
\end{proof}

\section{The inner form case}\label{sect:inner}
Let $D$ be the unique quaternion division algebra over $F$. Fix an element $\epsilon\in F^\times\setminus N_{E/F}E^\times$.
Let $G'=\GL_n(D)$ be an inner form of $\GL_{2n}(F)$ such that every $g\in G'$ in the image of $\GL_{2n}(E)$ is of the form $\begin{pmatrix}
A&\epsilon B\\ \bar{B}&\bar{A}
\end{pmatrix}$ where $A,B\in Mat_{n,n}(E)$. Let $\theta'$ be an involution of $G'$ given by
\[\theta'(g)=\begin{pmatrix}
\delta\\&-\delta
\end{pmatrix}g\begin{pmatrix}
\delta\\&-\delta
\end{pmatrix}^{-1} \]
for $g\in G'$. The fixed points coincide with $H_n=\GL_n(E)$.
\begin{thm}\label{thm:inner}
	Given any irreducible smooth admissible representation $\pi$ of $G'$, we have
	\[\dim\Hom_{H_n}(\pi,\mu)\leq1  \]
	for any character $\mu$ of $E^\times$.
\end{thm}
It has been proved by Guo in \cite{guo1997unique} when $\mu$ is trivial and $F$ is non-archimedean.
Similarly, we only need to prove the following
\begin{thm}
	Let $f$ be a bi-$(H_n,\mu)$-equivariant tempered generalized function on $G'$. Then
	\[f\Bigg(\begin{pmatrix}
	\epsilon\\&1
	\end{pmatrix}x^t\begin{pmatrix}
	\epsilon\\&1
	\end{pmatrix}^{-1} \Bigg)=f(x) \]
	for all $x\in G'$.\label{G'}
\end{thm}
The proof is very similar. The rest of this section is devoted to the proof of Theorem \ref{G'}. 
Let $\tilde{\mathcal{H}}_n'=\mathcal{H}_n\rtimes\langle\sigma'\rangle$ where $\sigma'$ acts on $\mathcal{H}_n=H_n\times H_n$ in the following way
\[\sigma'\cdot(g_1,g_2)=\Bigg(\begin{pmatrix}
\epsilon\\&1
\end{pmatrix}(g_2^{-1})^t\begin{pmatrix}
\epsilon\\&1
\end{pmatrix}^{-1}, \begin{pmatrix}
\epsilon\\&1
\end{pmatrix}(g_1^{-1})^t\begin{pmatrix}
\epsilon\\&1
\end{pmatrix}^{-1}\Bigg)=((g_2^{-1})^t,(g_1^{-1})^t) \]
for $g_i\in \GL_n(E)$.
Note that the map $g\mapsto \begin{pmatrix}
	\epsilon\\&1
\end{pmatrix}(g^{-1})^t\begin{pmatrix}
\epsilon\\&1
\end{pmatrix}^{-1}$ for $g\in\GL_n(D)$ is the analog of the Chevalley involution of $\GL_n(F)$, 
which takes an irreducible representation $\pi$ of $\GL_n(D)$ to its contragredient $\pi^\vee$ (see \cite[Theorem 3.1]{raghuram02} and \cite[Corollary 1.4]{adams}).
 Let $\mu\otimes\mu^{-1}$ be a character of $\mathcal{H}_n$ and $\tilde{\mu}$
be the character of $\tilde{\mathcal{H}}_n'$ twisted by the sign character as before. Let $\tilde{\mathcal{H}}_n'$ act on $G'$ by $(g_1,g_2)\cdot x=g_1xg_2^{-1}$ and
\[\sigma'\cdot x=\begin{pmatrix}
\epsilon\\&1
\end{pmatrix}x^t\begin{pmatrix}
\epsilon\\&1
\end{pmatrix}^{-1} \]
for $x\in G'$. Define 
\[L_n':=\{\begin{pmatrix}
0&\epsilon x\\\bar{x}&0
\end{pmatrix}|x\in Mat_{n,n}(E) \}\cong Mat_{n,n}(E) \]
and $\mathcal{N}_n'=\{\begin{pmatrix}
	0&\epsilon x\\ \bar{x}&0
\end{pmatrix}\in L_n'|\epsilon x\bar{x}\mbox{ is nilpotent} \}\cong\mathcal{N}_n$. Let $H_n$ acts on $L_n'$ by
\[g\cdot x=gx\bar{g}^{-1} \]
for $g\in H_n$. Then $\mathscr{C}_{\mathcal{N}_n'}(L_n')^{H_n,\mu_1}=0$ for any character $\mu_1$ of $E^\times/F^\times$ due to Proposition \ref{fourier}. 

Consider $x_k'=\begin{pmatrix}
&&\epsilon\mathbf{1}_k\\&\mathbf{1}_{n-k}\\
\mathbf{1}_k\\&&&\mathbf{1}_{n-k}
\end{pmatrix}$. Then $H_nx_k'H_n$ is a closed orbit in $H_n\backslash G'/H_n$. Thanks to \cite[Proposition 1.2]{guo1996cjm}, $H_nx'H_n$ exhaust all closed orbits in $H_n\backslash G'/H_n$ where
\[x'=\begin{pmatrix}
g_{1}&0&0&\epsilon g_{2}&0&0\\
0&0&0&0&\epsilon\mathbf{1}_k&0\\
0&0&\mathbf{1}_{n-k-\nu}&0&0&0\\
\bar{g}_{2}&0&0&\bar{g}_1&0&0\\
0&\mathbf{1}_k&0&0&0&0\\
0&0&0&0&0&\mathbf{1}_{n-k-\nu}\\
\end{pmatrix} \]
$g_1,g_2\in Mat_{\nu,\nu}(E)$ satisfying 
$$\begin{pmatrix}
	g_1&\epsilon g_2\\ \bar{g}_2&\bar{g}_1
\end{pmatrix}\begin{pmatrix}
\mathbf{1}_\nu\\&-\mathbf{1}_\nu
\end{pmatrix}\begin{pmatrix}
g_1&\epsilon g_2\\
\bar{g}_2&\bar{g}_1
\end{pmatrix}^{-1}\begin{pmatrix}
\mathbf{1}_\nu\\&-\mathbf{1}_\nu
\end{pmatrix}=\begin{pmatrix}
A&\epsilon B_A\\\bar{B}_A&{A}
\end{pmatrix} .$$ 
Here $A^2-\mathbf{1}_\nu=\epsilon B_A\bar{B}_A$ and $A\in Mat_{\nu,\nu}(F)$ is
semisimple without eigenvalues $\pm1$ and $AB_A=B_AA$.  Let $\tilde{\mathcal{H}}'_{n,x'}$ be the stablizer of $x'$.

\begin{proof}[Proof of Theorem \ref{G'}] Let $x',x_k'$, $\mathcal{H}_n'$, $\tilde{\mathcal{H}}_n'$ and $\tilde{\mu}$ be as above. It is equivalent to proving
	\[\mathscr{C}(G')^{\tilde{\mathcal{H}}_n',\tilde{\mu}}=0. \]
Applying Theorem \ref{thm2.1}, it  suffices to show that 
\begin{equation}\label{lem:G'}
\mathscr{C}(R(N_{H_nx'H_n}^{G'}))^{\tilde{\mathcal{H}}'_{n,x'},\tilde{\mu}}=0\Longrightarrow \mathscr{C}(Q(N_{H_nx' H_n}^{G'}))^{\tilde{\mathcal{H}}'_{n,x'},\tilde{\mu}}=0 \end{equation}
for all closed orbits $H_nx'H_n$. We separate them into two cases: $\nu=0$ and $k=0$.
\begin{itemize}
	\item Suppose that $\nu=0$. Then we have $N_{H_nx_k'H_n}^{G'}\cong Mat_{k,k}(E)\oplus L_{n-k}'$. The stabilizer of $x_k'$ in $\tilde{\mathcal{H}}_n'$ is $(\GL_{k}(E)\times \GL_{n-k}(E))\rtimes\langle\sigma'\rangle$, whose action on $N_{H_nx_k'H_n}^{G'}$
	is given by
	\[(g_1,g_2)\cdot(x,y)=(\bar{g}_1xg_1^{-1},g_2y\bar{g}_2^{-1}) \]
	and $\sigma'\cdot(x,y)=(x^t,\bar{y}^t)$
	for $g_1\in\GL_k(E)$, $g_2\in\GL_{n-k}(E),x\in Mat_{k,k}(E)$ and $y\in Mat_{n-k,n-k}(E)$. Then \eqref{lem:G'} holds (see the proof of Lemma \ref{nu=0}).
	\item Suppose that $k=0$. Suppose that $A$ is a scalar matrix and $A^2\neq\mathbf{1}_\nu$. Then
	$N_{H_nx'H_n}^{G'}\cong Mat_{\nu,\nu}(F)\oplus L_{n-\nu}'$. Thus
	 it is enough to prove the case for $\nu=n$ and $k=0$. 
	 Then
	\[\mathcal{H}_{n,x'}'\cong \GL_n(F)\rtimes\langle\sigma'\rangle \]
	$\tilde{\mu}|_{\mathcal{H}'_{n,x'}}=\chi$ is the sign character and $N_{H_nx'H_n}^{G'}\cong Mat_{n,n}(F)$.
	Here $\GL_n(F)$ acts on $Mat_{n,n}(F)$ by the inner conjugation.  Thus \eqref{lem:G'} holds (see the proof of Theorem \ref{thm:dist}).
\end{itemize}
This finishes the proof.
\end{proof}
Finally we can give a proof of Theorem \ref{thm:inner}.
\begin{proof}[Proof of Theorem \ref{thm:inner}]
	Thanks to Theorem \ref{G'} and \cite[Theorem 2.3]{sunzhu2011}, we have
	\[\dim\Hom_{\GL_{n}(E)}(\pi,\mu)\cdot\dim\Hom_{\GL_{n}(E)}(\pi^\vee,\mu^{-1})\leq1  \]
	for any irreducible smooth admissible representation $\pi$ of $\GL_{n}(D)$ and any  character $\mu$ of $E^\times$.
	Note that
	\[\dim\Hom_{\GL_{n}(E)}(\pi^\vee,\mu^{-1})=\dim \Hom_{\GL_{n}(E)}(\pi,\mu) \]
	since there exists an involution of $\GL_n(D)$  taking $\pi$ to its contragredient $\pi^\vee$  (see \cite[Theorem 3.1]{raghuram02} and \cite[Corollary 1.4]{adams}).
	Thus $\dim\Hom_{\GL_{n}(E)}(\pi,\mu)\leq 1 $, i.e. Theorem \ref{thm:inner} holds.
\end{proof}
\subsection*{Acknowlegement} The author would like to thank Dmitry Gourevitch for useful discussions. 
He also wants to thank Dipendra Prasad for explaining his work \cite{dipendraJNT}. 
He is also grateful for the anonymous referees for their careful
 reading of the manuscript and the numerous suggestions which greatly improve the exposition of this paper.
This work was partially supported by
the ERC, StG grant number 637912 and ISF grant 249/17. 

\bibliographystyle{amsplain}
\bibliography{galois}

\end{document}